# Uniformités et Continuity Spaces

Fleischer, Isidore et Giroux, Gaston

## Abstract

A semigroup A is an abelian semigroup with identity 0. A set of positives in A is an ordered down-directed set P containing with every r an element r/2 with r/2 + r/2 = r. A continuity space is an abstract set X equipped with a map d : XxX to A satisfying d(x, x) = 0 and

d(x, z) ≤d(x, y) + d(y, z). A quasi-uniform space is an abstract set X equipped with a filterbase of binary relations {U} such that each U contains the diagonal as well as $V \circ V$ for some V∈ {U}.
For each r∈ P, the set  $U(r) = \{(x,y) \mid d(x,y) < r \}$ is seen to be a quasi-uniform filterbase
on X .  Indeed, the down-directedness of P ensures that U(r) is a filterbase of oversets of the diagonal  and U(r) contains U(r/2)∘U(r/2).  One obtains a uniform filterbase by symmetrization,
i.e. by intersecting  the U(r) with the U(s) = {(y, x)|d(y, x) <s}.

## Résumé

Un semi-groupe A est un semigroupe abélien avec identité 0. Un ensemble d'éléments positifs est un ensemble ordonné dirigé vers le bas P qui contient avec chaque r un élément r/2 tel que r/2+r/2 = r. Un « continuity space » est un ensemble abstrait X muni d'une fonction d : XxX dans A satisfaisant d(x,x) = 0 et d(x, z) ≤d(x, y) + d(y, z). Un espace quasi-uniforme est un ensemble abstrait muni d'une base de filtre de relations binaires {U} tel que chaque U contient la diagonale ainsi que $V \circ V$   pour un  V∈ {U}. Pour chaque r∈ P, l'ensemble $U(r) = \{(x,y) \mid d(x,y) < r \}$ est une base de filtre d'une quasi-uniformité sur  X .  En effet le fait que P soit dirigé vers le bas nous assure que U(r) est une base de filtre dont les éléments contiennent la diagonale et qui est telle que U(r) contient U(r/2)∘U(r/2).  On obtient une base de filter d'une uniformité par symétrisassions,  i.e. en faisant les intersections  des U(r) avec les U(s) = {(y, x)|d(y, x) < s}.

Historiquement, les premières quasi-uniformités étudiées furent des uniformités; c'est-à-dire des quasi-uniformités possédant la propriété de symétrie :

(S)   $U \in \Psi \Rightarrow U^{-1} \in \Psi$

où   $U^{-1} = \{(x, y) \in XxX \mid (y, x) \in U \}$.

Par ailleurs un théorème célèbre, voir Weil (1937), identifie les espaces topologiques complètement réguliers aux espaces topologiques induits par les uniformités déterminées par certaines familles de fonctions continues sur X à valeurs réelles. Autrement dit les entourages sont de la forme $U_f = \{(x, y) : \mid f(x) - f(y) \mid < a \}$.

Rappelons que pour une uniformité $\Psi$ la topologie induite a pour système de voisinages les ensembles $U(x) = \{y \in X \mid (x, y) \in U \}$ où x varie dans X et U dans $\Psi$. Il revient au même, grâce à la propriété (S), de dire que le système de voisinages est $U^{-1}(y) = \{x \in X \mid (x, y) \in U \}$ où y varie dans X et U dans $\Psi$.

Mais il n'en est pas de même pour les quasi-uniformités ne satisfaisant pas la propriété (S). Elles apparaissent par paires et donnent lieu à 2 topologies.

Rappelons qu'une quasi-uniformité est un filtre de D(XxX) qui satisfait

(P)  $\forall U \in \Psi, \exists V \in \Psi, V \circ V \subseteq U$,

où D(XxX) est la famille des sous ensemble de XxX qui contiennent la diagonale.

Le concept d'une quasi-uniformité remonte à Czaszar (1960), d'après Murdeshwar-Naimpally (1966), et est plus fondamental que celui d'un espace topologique; on trouvera un bref aperçu historique dans Murdeshwar-Naimpally (1966). Une topologie n'est qu'un sous-produit d'une quasi-uniformité, même de plusieurs. De plus une quasi-uniformité induit un pré-ordre. C'est sur ce fait qu'est basée l'étude des espaces ordonnées uniformes de Nachbin (1965).

Dans son travail Kopperman (1988), l'auteur introduit la notion de métrique généralisée (qu'il aurait mieux value appeler une quasi-métrique généralisée). Bien entendu une telle métrique généralisée induit une quasi-uniformité par les entourages $U(r) = \{(x,y) \mid d(x,y) < r\}$. Cela résulte du fait que $U(r/2) \circ U(r/2) \subseteq U(r)$. En symétrisant, la quasi-uniformité devient une uniformité et le théoreme 11 de Kopperman (1988) un corollaire d'un résultat bien connu, voir Weil (1937).

Il est à noter que même si cela paraît étrange une topologie uniformisable induite par une uniformité, peut tout à la fois être induite par une quasi-uniformité. De même il n'y a pas correspondace 1-1 entre les « continuity spaces » et les quasi-uniformités.

## 6. Références

[1] Weil, A. (1937)Sur les espaces à structure uniforme et sur la topologie générale. Gauthier-Villars, Paris

[2] Czaszar, A. Fondements de la topologie générale Gauthier-Villars, Paris, 1960

[3] Murdeshwar, M. G. and Naimpally, S.A. (1966) Quasi-Uniform Topological Spaces. Preprints of Research Papers No4, Vol 2, University of Alberta QA 611.M87

[4] Nachbin, L. (1965) Topology and Order. D.Van Nostrand QA 611.N23

[5] Kopperman, R. (1988) All Topologies Come from Generalized Metrics. AMM, 95, 89-97

[6] Kelley, J. L. (1955) General Topology. van Nostrand